\documentclass[a4paper]{amsart}
\usepackage{amsmath,amsthm,amsfonts,amssymb,color}

\usepackage{a4wide}
\usepackage[dvips]{graphicx}
\usepackage{setspace}
\usepackage{url}
\theoremstyle{plain}
\newtheorem{thm}{Theorem}[section]
\newtheorem{prop}[thm]{Proposition}
\newtheorem{lem}[thm]{Lemma}
\newtheorem{defn}[thm]{Definition}
\newtheorem{cor}[thm]{Corollary}
\newtheorem{rem}[thm]{Remark}
\newtheorem{exam}[thm]{Example}

\newcommand{\R}{\mathbb{R}}

\newcommand{\bS}{\mathbb{S}}

\newcommand{\Z}{\mathbb{Z}}

\newcommand{\Pol}{\mathrm{Pol}}
\newcommand{\Harm}{\mathrm{Harm}}
\newcommand{\Span}{\mathrm{Span}}

\date{November 8, 2026}

\title{Null designs and switching in polynomial spaces}

\author[M.~Hirao]{Masatake Hirao}
\address[M.~Hirao]{Department of Information Science and Technology \\
    Tokyo University of Science \\
    2641 Yamazaki, Noda-shi, Chiba 278-8510, Japan}
\email{hirao@rs.tus.ac.jp}

\author[K.~Imamura]{Koji Imamura}
\address[K.~Imamura]{Institute of Mathematics for Industry (IMI) \\
    Kyushu University}
\email{k-imamura@imi.kyushu-u.ac.jp}

\author[H.~Nozaki]{Hiroshi Nozaki}
\address[H.~Nozaki]{Department of Mathematics Education \\
	Aichi University of Education \\
	1 Hirosawa, Igaya-cho, Kariya, Aichi 448-8542 \\
	Japan}
\email{hnozaki@auecc.aichi-edu.ac.jp}

\author[M.~Sawa]{Masanori Sawa}
\address[M.~Sawa]{Graduate School of System Informatics \\
    Kobe University \\
    1-1 Rokkodai, Nada, Kobe, Hyogo, 657-8501 \\
    Japan}
\email{sawa@people.kobe-u.ac.jp}

\author[T.~Tamaru]{Tomoki Tamaru}
\address[T.~Tamaru]{Graduate School of System Informatics \\
    Kobe University \\
    1-1 Rokkodai, Nada, Kobe, Hyogo, 657-8501 \\
    Japan}
\email{249x053x@stu.kobe-u.ac.jp}

\begin{document}

\begin{abstract}
Null designs are finite point configurations with real weights
whose polynomial moments vanish up to a given degree.
In Godsil's framework of polynomial spaces, we give elementary
proofs of sharp lower bounds on the number of points in the
Euclidean and spherical cases.
We characterize the spherical null designs attaining the bound
and give explicit formulas for their weights.
We then study how null designs can be used to replace part of
a weighted design while preserving its polynomial moments.
We give conditions ensuring nonnegative weights and apply
the minimum-size classification to switching on circles.
Examples on spheres and in Hamming and Johnson spaces
illustrate these results, including a replacement between
noncongruent spherical $5$-designs that preserves equal weights
and the number of points.
\end{abstract}

\maketitle

\section{Introduction}
\enlargethispage{3pt}
The eight vertices of a regular cube inscribed in $\bS^2$, each
with weight $1/8$, form a spherical $3$-design. Rotate the four
vertices of one face about its centre within its affine plane,
keeping the other four vertices fixed. For an angle not congruent
to zero modulo $\pi/2$, this replaces that face by four new points
and preserves the spherical $3$-design property. This is an instance
of the local deformation principle used in Sali's study of
rigidity \cite[Lemma~2.3]{Sali1994}; see Example~\ref{exam:cube}.

This face rotation preserves the spherical $3$-design property
because the original and rotated squares have the same polynomial
moments through degree $3$. More precisely, if $Y^+$ and $Y^-$
denote these squares, respectively, then
\[
 \sum_{y\in Y^+}f(y)-\sum_{y\in Y^-}f(y)=0
 \qquad\text{for every }f\in P_3(\R^3).
\]
Here $P_3(\R^3)$ denotes the space of polynomials of degree at most $3$.
More generally, a nonempty finite configuration on a sphere,
with nonzero real weights, is a spherical null $t$-design if its
weighted sum vanishes for every ambient polynomial of degree at
most $t$. Equivalently, the points with positive weights and those
with negative weights have the same weighted polynomial moments
through degree $t$, with the absolute values of the weights used
on the negative side.
Our aim is to study local configurations that can be exchanged
while preserving polynomial moments. We first establish universal
lower bounds on the support size of null designs and then classify
the spherical supports and weights attaining the bound.

The same principle of matching polynomial moments also appears
in combinatorial $t$-trades. The two collections have matching
incidence counts for subsets of
size at most $t$, or matching projections, with multiplicity,
onto at most $t$ coordinates. Consequently, their difference annihilates
the corresponding polynomial functions;
see, e.g.,~\cite{Billington2003,KhosrovshahiMaimaniTorabi1999}.
Combinatorial null designs and related objects have been studied by
Graver and Jurkat~\cite{GraverJurkat1973}, Graham et al.~\cite{GrahamLiLi1980},
and Deza et al.~\cite{DezaFranklSinghi1983}.
Frankl and Pach~\cite{FranklPach1983} proved the sharp lower bound
$2^{t+1}$ on the support size of a nontrivial combinatorial null
$t$-design. The corresponding bound for each sign class is $2^t$
\cite[Theorem~1]{DezaFranklSinghi1983}.
For general background on combinatorial and spherical designs, see
\cite{BethJungnickelLenz1999,BethJungnickelLenz1999_2,DelsarteGoethalsSeidel1977}.

To treat these moment identities in a common setting, we use Godsil's
polynomial spaces~\cite{Godsil1989,Godsil1993}, whose polynomial
function spaces are generated from a separation function.
Under a compatible spherical embedding, null $t$-designs in
a spherical polynomial space correspond exactly to spherical null
$t$-designs on the image. Thus geometric configurations and
embedded combinatorial trades can be studied as spherical null
designs and used in switching.
Moment-preserving replacement constructions in polynomial
spaces are also discussed by Nozaki and
Sawa~\cite[Section~3]{NozakiSawa2025}.
Support size means the number of points with nonzero weights.

Under the embedding assumptions in Section~\ref{sect:3}, we give
self-contained elementary proofs of the support-size bounds
$t+2$ in Euclidean polynomial spaces and $2(t+1)$ in spherical
polynomial spaces. In the spherical case, each of the positive and
negative sign classes has at least $t+1$ points.
Thus a moment-preserving complete replacement between disjoint
positively weighted parts requires at least $t+1$ points on each side.
The total bound applies to
spherical null designs on the standard sphere, in contrast with
the exponential bound in the combinatorial setting above.
The bounds on the total support size and, for $t\ge2$, the
characterization of minimum-size spherical supports as sets of
$2t+2$ distinct points on a circle follow from Couvreur's results
\cite{Couvreur2012}; their relation to our setting is explained
in Section~\ref{sect:3}.
In a general polynomial space, every sufficiently large finite
subset contains the support of a null design.

The separate lower bounds on the numbers of positively and
negatively weighted points in a spherical null design also follow
from the uniqueness result for nonnegative measures with small
support on the sphere \cite[Section~5]{BendoryDekelFeuer2015}.
The spherical bound gives an affirmative answer to the conjecture
stated by Bannai~\cite{Bannai2013};
see also Campena~\cite{Campena2016} for subsequent work.

For spherical null designs, we next determine the supports and
weights attaining the bound $2(t+1)$, for $t\ge1$
(Theorem~\ref{thm:minimal-spherical-null-classification}).
Using classical trigonometric interpolation, we give explicit
formulas for their weights.
The supports are precisely the sets of $2(t+1)$ distinct points
on a circle contained in the sphere.
Their null weights are unique up to
a nonzero scalar and alternate in sign in cyclic order.
We also characterize the case of equal absolute weights as
the difference of two disjoint regular $(t+1)$-gons
on the same circle, with arbitrary relative rotation subject
to disjointness (Corollary~\ref{cor:minimal-null-equal-weights}).
The two squares in the opening example form an eight-point null
$3$-design of affine dimension $2$, attaining the bound.
Together, the elementary proofs and explicit weight formulas give
a self-contained account of minimum-size spherical null designs.
In the equal-absolute-weight case, the classification recovers the
regular-polygon structure underlying the opening rotation construction.

To use these null configurations for switching, their placement
and weights must be compatible with the design to be modified.
We formulate these conditions for weighted designs in spherical
polynomial spaces and apply the explicit null weights to switching
on circles (Corollary~\ref{cor:circle-switching}).
The discussion following the corollary distinguishes partial
modification from complete replacement. Trades from Hamming and Johnson spaces
provide another source of configurations through suitable spherical
embeddings. In both cases, the null-design condition ensures that
polynomial moments are preserved.

The paper is organized as follows.
In Section~\ref{sect:2}, we recall the basic definitions of polynomial
spaces and designs and introduce null designs.
In Section~\ref{sect:3}, we discuss support-size bounds and
existence results, including the criterion for a prescribed support,
explain their relation to earlier work, and describe the
minimum-size spherical null designs.
In Section~\ref{sect:4}, we develop the switching construction for
spherical polynomial spaces.
In Section~\ref{sect:5}, we realize Hamming and Johnson trades
as spherical null designs and illustrate their use in switching.
An equal-weight replacement in a $240$-point spherical $5$-design
yields a noncongruent design (Example~\ref{exam:240point-switching}).
Finally, in Section~\ref{sect:6}, we give further remarks,
including a brief discussion of the connection with the
Prouhet--Tarry--Escott problem in number theory.

\section{Preliminaries}
\label{sect:2}

In this section, we review basic definitions and results on
polynomial spaces and designs.
We follow Godsil~\cite{Godsil1989,Godsil1993} and Section~2 of Nozaki and Sawa~\cite{NozakiSawa2025}.

\subsection{Polynomial space}
Let $\Omega$ be a possibly infinite set.
A \textit{separation function} $\rho$ on $\Omega$
is a function from $\Omega \times \Omega$ to $\R$.
For each $a\in\Omega$, define $\rho_a:\Omega\to\R$ by
\[
\rho_a(\xi)=\rho(a,\xi).
\]
For $f \in \mathbb{R}[x]$,
the \textit{zonal polynomial} with respect to $a \in \Omega$
is $f \circ \rho_a: \Omega \to \R$.
Let $Z(\Omega, r)$ denote the linear space spanned by all zonal polynomials of degree at most $r$:
\[
Z(\Omega, r) := \Span_{\R} \{ f \circ \rho_a \mid a \in \Omega, \; f \in \R[x], \; \deg f \leq r \}.
\]
Let $\Pol(\Omega,0)$ be the space of constant functions on $\Omega$, which coincides with $Z(\Omega,0)$, and set
$\Pol(\Omega, 1) = Z(\Omega, 1)$.
We define inductively
\[
\Pol(\Omega, r + 1) = \Span_{\R} \{ gh \mid g \in \Pol(\Omega, r), \; h \in \Pol(\Omega, 1)  \}.
\]
While $Z(\Omega,r)$ is spanned by zonal polynomials of degree at most $r$, the space $\Pol(\Omega,r)$ also contains products of zonal polynomials whose total degree is at most $r$.

Let $\Pol(\Omega) = \bigcup_{r \ge 0} \Pol(\Omega, r)$.
An element of $\Pol(\Omega)$ is called a \textit{polynomial} (function) on $\Omega$.
An element of $\Pol(\Omega, r) \setminus \Pol(\Omega, r-1)$
is called a \textit{polynomial of degree $r$} on $\Omega$.
The \textit{dimension} of a polynomial space is $\dim \Pol(\Omega, 1) - 1$.

The space $\Pol(\Omega, r)$ depends on the separation function $\rho$.
To emphasize the dependence on the separation function $\rho$, we write $\Pol_{\rho}(\Omega, r)$.

The concept of polynomial space was defined by Godsil~\cite{Godsil1989,Godsil1993}.
\begin{defn}[Polynomial space]
Let $\rho$ be a separation function on a set $\Omega$.
Let $\langle f, g \rangle$ be an inner product on $\Pol(\Omega)$.
The pair $(\Omega, \rho)$ is called a \textit{polynomial space}
if the following conditions hold.
\begin{enumerate}
\item For any $x, y \in \Omega$, $\rho(x, y) = \rho(y,x)$.
\item The dimension of $\Pol(\Omega, r)$ is finite for any $r$.
\item For any $f, g \in \Pol(\Omega)$, $\langle f, g \rangle = \langle 1, fg \rangle$.
\item For any $f\in\Pol(\Omega)$ satisfying $f(x)\ge0$
for every $x\in\Omega$, we have $\langle 1,f\rangle\ge0$,
with equality if and only if $f=0$.
\end{enumerate}
\end{defn}

We normalize the inner product $\langle \cdot, \cdot \rangle$
so that $\langle 1, 1 \rangle = 1$.
We usually use the following inner product on $\Pol(\Omega)$
for a finite underlying set $\Omega$:
\[
\langle f, g \rangle = \sum_{x \in \Omega} f(x) g(x) \mu (x)
\]
for $f, g \in \Pol (\Omega)$,
where $\mu$ is a probability measure on $\Omega$
such that $\mu(x)>0$ for every $x\in\Omega$.

An affine transformation with nonzero slope of the separation function
does not change the associated polynomial spaces.
More precisely, two separation functions $\rho$ and $\sigma$ on $\Omega$ are said to be \textit{affinely equivalent} if there exist
$a\in\R_{\ne 0} :=\R\setminus \{0\}$ and $b\in\R$ such that
\[
\rho(x,y)=a\sigma(x,y)+b
\]
for all $x,y\in\Omega$.
If $\rho$ and $\sigma$ are affinely equivalent, then
\[
\Pol_{\rho}(\Omega,r)=\Pol_{\sigma}(\Omega,r)
\]
for every $r\geq0$.
Thus, the polynomial spaces depend only on the affine equivalence class
of the separation function.

We next recall the definitions of Euclidean and spherical polynomial
spaces. In these spaces, polynomial functions are restrictions of
ordinary polynomials on Euclidean space. A spherical polynomial space
is a Euclidean polynomial space admitting an embedding into a unit sphere.

\begin{defn}[Euclidean and spherical polynomial spaces]
A polynomial space $(\Omega, \rho)$ is called
\textit{Euclidean} (resp. \textit{spherical}) if there exists an injection $\tau$ from $\Omega$ to the Euclidean space $\R^d$ (resp. the unit sphere $\bS^{d-1}$)
for some $d$ and a separation function $\rho'$ affinely equivalent to $\rho$, such that for any $x, y \in \Omega$,
\[
\rho'(x, y) = (\tau(x), \tau(y)),
\]
where $(\cdot, \cdot)$ is the usual inner product of $\R^d$.
\end{defn}

The spaces $Z(\Omega,r)$ and $\Pol(\Omega,r)$ need not coincide
in a general polynomial space. For Euclidean polynomial spaces,
Godsil proved the following theorem.

\begin{thm}[Zonal polynomial representation, \cite{Godsil1989,Godsil1993}]
A Euclidean polynomial space $(\Omega, \rho)$
satisfies $Z(\Omega, r) = \Pol(\Omega, r)$
for any $r \ge 0$.
\end{thm}

The preceding theorem describes $\Pol(\Omega,r)$ in terms of
the separation function. The following theorem expresses
$\Pol(\Omega,r)$ as restrictions of ordinary polynomials under
an embedding that realizes this function as the Euclidean inner product.

We write $\R_r[x_1,\ldots,x_n]$, or equivalently $P_r(\R^n)$,
for the space of real polynomials in $n$ variables
of total degree at most $r$.
\begin{thm}[Polynomial representation under an embedding, \cite{Godsil1989,Godsil1993}]
\label{thm:Godsil}
Let $(\Omega, \rho)$ be a Euclidean polynomial space
of dimension $n$.
Let $\tau$ be an injection from $\Omega$ to $\R^n$
such that for any $x, y \in \Omega$, $\rho(x, y) = (\tau(x), \tau(y))$.
Then for each $r \ge 0$,
\[
\Pol (\Omega, r) = \{ f \circ \tau: \Omega \to \R \mid
f \in \R_r[x_1, \ldots, x_n]\}.
\]
\end{thm}

For a polynomial space $(\Omega, \rho)$,
an injection $\tau: \Omega \to \R^n$ is said to be
\textit{affinely compatible} with $\rho$
if there exist $a \in \R_{\neq 0}$ and $b \in \R$ such that
for any $x, y \in \Omega$, $ \rho(x, y) = a(\tau(x), \tau(y)) + b$.
If there exists an affinely compatible injection $\tau$,
then separation functions $\rho$ and
\[
\sigma(x, y) = \frac{1}{a} \rho (x, y) - \frac{b}{a}
= (\tau(x), \tau(y))
\]
are affinely equivalent,
and hence $\Pol_{\rho}(\Omega, r) = \Pol_{\sigma}(\Omega, r)$.

By this observation and the preceding theorem, we obtain the
following corollary.

\begin{cor}[cf.~\cite{Godsil1993,NozakiSawa2025}]
\label{cor:polynomial-representation}
Let $(\Omega,\rho)$ be a Euclidean (resp.\ spherical) polynomial space
of dimension $n$.
Suppose that there exists an injection
\[
    \tau:\Omega\to\R^n
    \quad
    (\text{resp. } \tau:\Omega\to\bS^{n-1})
\]
that is affinely compatible with $\rho$.
Then, for each $r\geq 0$,
\[
    \Pol(\Omega,r)
    =
    \left\{
        f\circ\tau:\Omega\to\R
        \;\middle|\;
        f\in\R_r[x_1,\ldots,x_n]
    \right\}.
\]
\end{cor}

\subsection{Weighted design and null design}

We now define weighted designs and null designs in a polynomial space.
\begin{defn}[Weighted $t$-design, \cite{NozakiSawa2025}]
Let $(\Omega, \rho)$ be a polynomial space.
Let $X$ be a finite subset of $\Omega$
and $w: X \to \R_{> 0}$ be a positive weight function.
The pair $(X, w)$ is called a \textit{weighted $t$-design} if
\[
\sum_{x \in X} w(x) f(x) = \langle 1, f \rangle
\quad\text{for all } f\in\Pol_{\rho}(\Omega,t).
\]
\end{defn}

\begin{defn}[Null $t$-design]
Let $(\Omega, \rho)$ be a polynomial space.
Let $X$ be a nonempty finite subset of $\Omega$
and $w: X \to \R_{\ne 0}$ be a weight function.
The pair $(X, w)$ is called a \textit{null $t$-design in $(\Omega, \rho)$} if
\[
\sum_{x \in X} w(x) f(x) = 0
\quad\text{for all } f\in\Pol_{\rho}(\Omega,t).
\]
\end{defn}
By definition, all null designs considered here have nonempty
support and nonzero weights at every support point.
A null design has \textit{strength exactly $t$} if it is a null
$t$-design but not a null $(t+1)$-design.
Taking $f=1$ in the defining condition gives $\sum_{x\in X}w(x)=0$.
Thus both signs occur, and the total positive weight equals the
sum of the absolute values of the negative weights.

For a nonempty finite set $X\subset\R^n$ and a weight function
$w:X\to\R_{\ne0}$, we call $(X,w)$ a \textit{Euclidean null $t$-design} if
\[
 \sum_{x\in X}w(x)f(x)=0
 \qquad\text{for all }f\in P_t(\R^n).
\]
When $X\subset\bS^{n-1}$, we call it a \textit{spherical null $t$-design}.

The difference of two distinct positive weighted $t$-designs in
the same polynomial space is a null $t$-design: extend their weights
by zero to the union of their supports, subtract the weights,
and remove points with zero coefficients. Conversely, adding a null $t$-design
to a weighted $t$-design preserves its prescribed moments;
if the resulting weights are nonnegative, deleting zero weights
yields a positive weighted $t$-design.

For polynomial spaces admitting an affinely compatible embedding,
the weighted and null design conditions can be expressed in terms of ordinary
polynomials.  In particular, the following characterization applies
to both Euclidean and spherical polynomial spaces.

\begin{prop}
\label{prop:design-polynomial-space}
Under the assumptions of
Corollary~\ref{cor:polynomial-representation},
let $X$ be a nonempty finite subset of $\Omega$.
Then the following hold.
\begin{enumerate}
\item For a positive weight function $w:X \to \R_{>0}$, $(X,w)$ is a weighted $t$-design if and only if
\[
\sum_{x\in X} w(x) f(\tau(x)) = \left\langle 1, f\circ\tau\right\rangle \qquad \text{for all }f\in\R_t[x_1,\ldots,x_n].
\]

\item For a weight function $w:X \to \R_{\ne 0}$, $(X,w)$ is a null $t$-design in $(\Omega,\rho)$ if and only if
\[
\sum_{x\in X}w(x)f(\tau(x))=0 \qquad \text{for all }f\in\R_t[x_1,\ldots,x_n].
\]
\end{enumerate}
\end{prop}

\begin{proof}
By Corollary~\ref{cor:polynomial-representation},
\[
    \Pol_\rho(\Omega,t)
    =
    \{f\circ\tau:\Omega\to\R
      \mid f\in\R_t[x_1,\ldots,x_n]\}.
\]
Hence the assertion follows immediately from the definitions of
weighted $t$-designs and null $t$-designs.
\end{proof}

Proposition~\ref{prop:design-polynomial-space} identifies a null
$t$-design in a spherical polynomial space with a spherical null
$t$-design on its image under a compatible embedding $\tau$,
with weights $\widehat w(\tau(x))=w(x)$. We use this identification
throughout, retaining the original space when its combinatorial
structure is relevant. The null condition involves only polynomial
evaluations. For a positive weighted design, the image instead
reproduces the averages $\langle1,f\circ\tau\rangle$ from the
original space. These need not be the averages of $f$ over the
sphere. The image is a weighted spherical $t$-design precisely
when the two averages agree for every polynomial of degree at
most $t$.

\section{Bounds, classification, and existence of null designs}
\label{sect:3}

\subsection{Lower bounds for null designs}

We establish general lower bounds on the support size of null
$t$-designs in Euclidean and spherical polynomial spaces.
Our argument uses suitable polynomials constructed from the
separation function.

Let $(X, w)$ be a null $t$-design in $(\Omega, \rho)$.
Partition $X$ according to the sign of the weight function $w$, i.e.,
\[
X^+ := \{ x \in X \mid w(x) > 0 \}, \qquad
X^- := \{ x \in X \mid w(x) < 0 \}.
\]

\begin{thm}[Sign-class bounds]
\label{thm:main1}
\begin{enumerate}
\item
Let $(\Omega, \rho)$ be a spherical polynomial space of dimension $n$.
Suppose there exists an injection $\tau: \Omega \to \bS^{n-1}$ that is affinely compatible with $\rho$.
Let $(X, w)$ be a null $t$-design in $(\Omega, \rho)$.
Then
\[
\min \{ |X^+|, |X^-| \} \geq t + 1, \; \quad
\textrm{and} \; \quad |X| \geq 2 t + 2.
\]

\item
Let $(\Omega, \rho)$ be a Euclidean polynomial space of dimension $n$.
Suppose there exists an injection $\tau: \Omega \to \R^n$ that is affinely compatible with $\rho$.
Let $(X, w)$ be a null $t$-design in $(\Omega, \rho)$.
Then
\[
\min \{ |X^+|, |X^-| \} \geq \lfloor t/2 \rfloor + 1, \; \quad
\textrm{and} \; \quad  |X| \geq 2 \lfloor t/2 \rfloor + 2.
\]
\end{enumerate}
\end{thm}

\begin{proof}
We prove (1) first.
Define $\rho'_x(y) := (\tau(x), \tau(y))$.
Then, by the affine compatibility of $\tau$ with $\rho$,
$\rho'$ is affinely equivalent to $\rho$.
Moreover, $\Pol_{\rho}(\Omega, t) = \Pol_{\rho'}(\Omega, t)$.

We first show that $|X^+|\ge t+1$.
Suppose, to the contrary, that $|X^+| \le t$.
Consider the polynomial function
\begin{equation*}
p_+(y) := \prod_{x \in X^+}
\left (2 - 2 \rho'_x(y) \right ) \in \Pol_{\rho}(\Omega, t).
\end{equation*}
Its membership in $\Pol_\rho(\Omega,t)$ follows from
Corollary~\ref{cor:polynomial-representation}.
Since $p_+(y) = 0$ for $y \in X^+$ and $p_+(y) > 0$ for $y \in X^-$, while $w(y) < 0$ for $y \in X^-$, we obtain
\[
\sum_{y\in X}w(y)p_+(y)
=
\sum_{y\in X^{-}}w(y)p_+(y)<0.
\]

On the other hand, since
$\deg p_+ \le t$,
the null design condition yields
\[
\sum_{y\in X}w(y)p_+(y)=0,
\]
a contradiction.
Hence
\[
|X^+|\ge t+1.
\]

Applying the same argument to
$p_-(y) := \prod_{x \in X^-} \left (
1 - \rho_x'(y) \right )$
shows that $|X^-|\ge t+1$.
Finally, we also obtain $|X| \geq 2 \min \{ |X^+|, |X^-| \} \ge 2t + 2$.

\bigskip

For (2), we use the same argument with a different choice of
$p_+(y)$.
Suppose, to the contrary, that $|X^+| \leq \lfloor t/2 \rfloor$.
Consider the polynomial function
\begin{align*}
p_+ (y) &:= \prod_{x \in X^+} \| \tau(x) - \tau(y) \|^2 \\
&= \prod_{x \in X^+}
\left \{ (\tau(x), \tau(x)) - 2 (\tau(x), \tau(y)) +
(\tau(y), \tau(y)) \right \} \\
&= \prod_{x \in X^+}
\left \{ \rho'_x(x) - 2 \rho'_x(y) + \rho'_y(y) \right \}
\in \Pol_{\rho}(\Omega, 2 \lfloor t/2 \rfloor)
\subset \Pol_{\rho}(\Omega, t).
\end{align*}
The membership assertions follow from
Corollary~\ref{cor:polynomial-representation}.
The rest of the proof is the same as in (1).
\end{proof}

\begin{rem}
On the unit sphere, $2-2(u,v)=\|u-v\|^2$, so the two parts of
Theorem~\ref{thm:main1} use the same product of squared distances.
In part~(1), $p_+$ is the pullback under $\tau$ of this product
on the sphere, with factors indexed by $\tau(X^+)$.

The same polynomial construction shows that a finite nonnegative
measure supported on at most $t$ points of the sphere is uniquely
determined, among finite nonnegative measures on the sphere, by
its polynomial moments through degree $t$.
This uniqueness result is also established in
\cite[Section~5]{BendoryDekelFeuer2015}; applied to the positive
and negative parts of a null design, it gives the spherical
sign-class bounds in Theorem~\ref{thm:main1}(1).
\end{rem}

By Theorem~\ref{thm:main1}, we obtain the following corollary,
which gives an affirmative answer to the conjecture
stated by Bannai~\cite{Bannai2013}.
As discussed below, the same bound also follows from
Couvreur's general results~\cite{Couvreur2012}.

\begin{cor}
\label{cor:Bannai-conjecture}
Let $(X,w)$ be a spherical null $t$-design on
$\bS^{n-1}$, where $n\geq 2$. Then
\[
    |X| \geq 2(t+1).
\]
\end{cor}

We next improve the lower bound on $|X|$ in
Theorem~\ref{thm:main1}(2) when $t$ is odd, by using evaluation functionals.
For each point $x\in X$,
define the evaluation functional
\[
L_x:\Pol(\Omega, t) \to\mathbb R,
\quad
L_x(f)=f(x).
\]

Thus every point determines a linear functional on the polynomial space $\Pol(\Omega, t)$.
The family $\{L_x:x\in X\}$
determines all evaluations of polynomials on the finite set $X$.

We note that
$(X,w)$ is a null $t$-design in $(\Omega, \rho)$
if and only if
\[
\sum_{x\in X}w(x)L_x=0
\]
holds in the dual space $\Pol(\Omega, t)^*$.
Hence a nonempty finite set $X\subset\Omega$ contains the support
of a null $t$-design if and only if the family $(L_x)_{x\in X}$
is linearly dependent.
Indeed, the support of any nonzero linear dependence gives
such a null design.
Consequently, whenever null $t$-designs exist, their minimum
support size equals the smallest number of points whose
associated evaluation functionals are linearly dependent.

\begin{thm}[Euclidean support bound]
Let $(\Omega, \rho)$ be a Euclidean polynomial space of dimension $n$.
Suppose there exists an injection $\tau: \Omega \to \R^n$ that is affinely compatible with $\rho$.
Let $(X, w)$ be a null $t$-design in $(\Omega, \rho)$.
Then,
\[
|X| \geq t + 2.
\]
\end{thm}

\begin{proof}
Let $\{ x_1, \ldots, x_m \}$ be a subset of
$\Omega$ with $m \leq t + 1$.
Set $\rho'(x,y)=(\tau(x),\tau(y))$.
Then $\rho'$ is affinely equivalent to $\rho$.
For each $i \neq j$, define
\[
\ell_{ij}: \Omega \to \R, \quad
\ell_{ij}(x) = (\tau(x), \tau(x_{i}) - \tau(x_j))
= \rho'_{x_i}(x) - \rho'_{x_j}(x).
\]
Since $\tau$ is injective, for $i\ne j$ we have
\[
\ell_{ij}(x_i)-\ell_{ij}(x_j)
=\|\tau(x_i)-\tau(x_j)\|^2>0.
\]
For each $i = 1, \ldots, m$, let $p_i$
be the interpolation function defined by
\[
p_i(x) = \prod_{\substack{1 \leq j \leq m \\ j \neq i}}
\frac{\ell_{ij}(x) - \ell_{ij}(x_j)}{\ell_{ij}(x_i) - \ell_{ij}(x_j)} \in \Pol_{\rho}(\Omega, t).
\]

Suppose $\sum_{i=1}^m c_iL_{x_i}=0$ in the dual space $\Pol(\Omega, t)^*$.
Evaluating both sides at $p_i$ yields
\[
0 = \sum_{j=1}^m c_jp_i(x_j) = c_i, \quad i = 1, \ldots, m.
\]
Hence every coefficient vanishes, so the evaluation functionals are linearly independent.
Thus no null $t$-design can have at most $t+1$ support points, and hence $|X|\ge t+2$.
\end{proof}

\begin{rem}
To relate our setting to Couvreur's results, use homogeneous
coordinates $[Z_0:Z_1:\cdots:Z_n]$ on real projective space
$\mathbb P^n(\R)$, and put $Z=(Z_1,\ldots,Z_n)$.
The map $a\mapsto[1:a]$ identifies $\R^n$ with the affine chart
$Z_0\ne0$. Under Corollary~\ref{cor:polynomial-representation},
homogenization gives
\[
 \widetilde f(Z_0,Z)=Z_0^t f(Z/Z_0),\qquad
 \widetilde f(1,\tau(x))=f(\tau(x))
\]
for $f\in P_t(\R^n)$, where the first expression is expanded
as a homogeneous polynomial of degree $t$, also defined at $Z_0=0$.
Thus null relations correspond to degree-$t$ evaluation relations
at $[1:\tau(x)]$. Couvreur's Proposition~6.1 gives the bound $t+2$;
for $t\ge2$, Proposition~7.3 gives the spherical bound $2t+2$,
since no three distinct points on a sphere are collinear
\cite{Couvreur2012}. The cases $t=0,1$ follow directly from
constant and affine linear polynomials. These applications over
$\R$ are compatible with \cite[Lemma~2.13]{Couvreur2012}.

For $t\ge2$, a spherical null support of size $2t+2$ has
minimally dependent evaluations by the spherical bound.
Proposition~7.4 of \cite{Couvreur2012} then implies coplanarity,
since $2t+2\le3t$, and hence containment in a circle.
This gives the conclusion of Proposition~\ref{thm:spherical-equality}
below; for $t=1$, it follows directly from affine dependence.
Conversely, for $t\ge2$, Couvreur's Theorem~5.2 and Lemma~7.1 show that
$2t+2$ distinct points on a circle have minimally dependent
degree-$t$ evaluations.
Our proofs establish these conclusions directly, without using
projective geometry, and also give separate bounds for the two
sign classes. We then determine the null weights explicitly and
use the resulting configurations in switching constructions.
\end{rem}

\subsection{Minimum-size spherical null designs}
\label{subsect:minimum-spherical}

We now consider the equality case in the spherical bound of
Theorem~\ref{thm:main1}.
We show that, for $t\ge1$, every spherical null $t$-design
$(X,w)$ with $|X|=2(t+1)$ lies on a circle,
not necessarily a great circle.
For $t=0$, the minimum-size null designs consist of two distinct
points with opposite nonzero weights. We therefore restrict the
following classification to $t\ge1$.
For $X\subset\R^n$, let $\operatorname{aff}X$ denote the smallest
affine subspace containing $X$, and write
$\dim\operatorname{aff}X$ for its dimension.

\begin{lem}
\label{lem:sphere}
Any affine line intersects $\bS^{n-1}$ in at most two points.
In particular, no three distinct points on $\bS^{n-1}$ are collinear.
\end{lem}

\begin{proof}
Write an affine line as $\ell(s)=as+b$, where $a\ne0$.
We have $\ell(s)\in\bS^{n-1}$ if and only if
\[
 \|a\|^2s^2+2(a,b) s+\|b\|^2-1=0.
\]
Since $\|a\|^2>0$, this equation has at most two real roots.
\end{proof}

\begin{lem}
\label{lem:separating-polynomial}
Suppose that $X\subset\bS^{n-1}$ satisfies
\[
 |X|=2(t+1),\qquad \dim\operatorname{aff}X\ge3.
\]
Then, for each $p\in X$, there exists a polynomial
$F_p\in P_t(\R^n)$ such that
\[
 F_p(p)>0,\qquad F_p(x)=0\quad(x\in X\setminus\{p\}).
\]
\end{lem}

\begin{proof}
Fix $p\in X$.
Since $|X|\ge4$, we may choose two distinct points
$a,b\in X\setminus\{p\}$.
By Lemma~\ref{lem:sphere}, the points $p,a,b$ are not collinear,
so $A:=\operatorname{aff}\{p,a,b\}$ is a two-dimensional
affine subspace.
Since $\dim\operatorname{aff}X\ge3$, we may choose
$c\in X\setminus A$.
Set $G_1:=\{a,b,c\}$ and $B:=\operatorname{aff}G_1$.
Both $A$ and $B$ are two-dimensional affine subspaces, and
$\operatorname{aff}\{a,b\}\subset A\cap B$.
Moreover, $A\ne B$ because $c\notin A$.
Thus $A\cap B=\operatorname{aff}\{a,b\}$, and hence
$p\notin B=\operatorname{aff}G_1$.

If $t\ge2$, partition the remaining
$|X\setminus\{p,a,b,c\}|=2(t-1)$ points into disjoint pairs:
\[
 X\setminus\{p,a,b,c\}=G_2\sqcup\cdots\sqcup G_t,
 \qquad |G_i|=2\quad(2\le i\le t).
\]
By Lemma~\ref{lem:sphere}, we have
$p\notin\operatorname{aff}G_i$ for each $i=2,\ldots,t$.
For $t=1$, there are no remaining points, and only $G_1$ is needed.

Now let $\Delta=d_0+W$ be an affine subspace not containing $p$,
and let $v$ be the orthogonal projection of $p-d_0$ onto $W^\perp$.
Then $v\ne0$, and the affine linear function
$\alpha(x):= (v,x-d_0)$
satisfies $\alpha(p)=\|v\|^2>0$ and $\alpha|_\Delta=0$.
Applying this construction to $\Delta=\operatorname{aff}G_i$,
we obtain, for each $i=1,\ldots,t$, an affine linear function
$\alpha_i$ satisfying
\[
 \alpha_i(p)>0,\qquad \alpha_i(g)=0\quad(g\in G_i).
\]
Define the polynomial
\[
    F_p(x):=\prod_{i=1}^t\alpha_i(x),
\]
which has degree at most $t$.
Every $x\in X\setminus\{p\}$ belongs to exactly one $G_i$,
so $\alpha_i(x)=0$ and hence $F_p(x)=0$.
Finally, since $\alpha_i(p)>0$ for every $i$, we have $F_p(p)>0$.
\end{proof}

\begin{prop}
\label{thm:spherical-equality}
Let $(X,w)$ be a spherical null $t$-design on $\bS^{n-1}$
with $t\ge1$, and suppose that $|X|=2(t+1)$.
Then
\[
 \dim\operatorname{aff}X=2.
\]
In particular, $X$ lies on a circle obtained as the intersection
of $\bS^{n-1}$ with a two-dimensional affine subspace.
\end{prop}

\begin{proof}
Suppose that $\dim\operatorname{aff}X\ge3$.
Choose $p\in X$, and let $F_p$ be the polynomial given by
Lemma~\ref{lem:separating-polynomial}.
Since $\deg F_p\le t$, the null design condition gives
\[
 0=\sum_{x\in X}w(x)F_p(x)=w(p)F_p(p)\ne0,
\]
a contradiction.
Here $w(p)\ne0$ by the definition of a null design.
Thus $\dim\operatorname{aff}X\le2$.
On the other hand, $t\ge1$ implies $|X|=2(t+1)\ge4$,
so Lemma~\ref{lem:sphere} yields
$\dim\operatorname{aff}X\ge2$.
Therefore $\dim\operatorname{aff}X=2$.

To prove the last assertion, put $\Pi:=\operatorname{aff}X$.
Let $c$ be the point of $\Pi$ closest to the origin.
Then $\Pi=c+U$ for a two-dimensional linear subspace
$U\subset\R^n$ orthogonal to $c$.
For $x=c+u\in\Pi$, we have $c\perp u$, and hence
\[
 \|x\|^2=\|c\|^2+\|u\|^2.
\]
Therefore,
\[
 \Pi\cap\bS^{n-1}
 =\{c+u\mid u\in U,\ \|u\|^2=1-\|c\|^2\}.
\]
Since $X\subset\Pi\cap\bS^{n-1}$ and $|X|\ge4$,
we have $1-\|c\|^2>0$.
Setting $r:=\sqrt{1-\|c\|^2}$, we obtain
\[
 \Pi\cap\bS^{n-1}
 =\{c+u\mid u\in U,\ \|u\|=r\},
\]
which is a circle in $\Pi$ containing $X$.
\end{proof}

By the preceding proposition, it remains to consider points on a circle.
The converse is obtained by classical trigonometric interpolation.
For the corresponding trigonometric Vandermonde determinant, see
\cite[Lemma~2.3]{Karp2019}. We include a proof to describe the weights explicitly.

\begin{thm}[Classification of minimum-size spherical null designs]
\label{thm:minimal-spherical-null-classification}
Let $n\geq 2$ and $t\geq 1$, and let
$X=\{x_1,\ldots,x_{2t+2}\}\subset\bS^{n-1}$ consist of distinct points.
Then $X$ supports a spherical null $t$-design if and only if $X$ is
contained in a circle obtained by intersecting $\bS^{n-1}$ with a
two-dimensional affine subspace.

In this case, the circle is unique, and the null weight function
is unique up to a nonzero real constant multiple.
With the points indexed in cyclic order around the circle, all null weight functions are given by
\begin{equation}
\label{eq:minimal-null-distance-weights}
 w(x_j)=\alpha\,\frac{(-1)^{j-1}}
 {\displaystyle\prod_{\ell\ne j}\|x_j-x_\ell\|},
 \qquad \alpha\in\R_{\ne0}.
\end{equation}
Here the norm is the ambient Euclidean norm.
Their signs alternate in cyclic order. In particular,
$|X^+|=|X^-|=t+1$, and every such null design has strength exactly $t$.
\end{thm}

\begin{proof}
The only-if part follows from Proposition~\ref{thm:spherical-equality}. Any three
distinct points on the sphere are noncollinear, so a circle containing
$X$ must lie in the plane $\operatorname{aff}X$. Thus the circle is unique.

Conversely, suppose that $X$ lies on such a circle. We write
\[
 \gamma(\theta)=c+r(\cos\theta\,u+\sin\theta\,v),
\]
where $u,v$ are orthonormal, $c\perp u,v$, $r>0$, and
$\|c\|^2+r^2=1$. Choose the angular origin and orientation so that
\[
 x_j=\gamma(\theta_j),\qquad
 0\leq\theta_1<\cdots<\theta_{2t+2}<2\pi.
\]
The restrictions to this circle of polynomials in
$P_t(\R^n)$, expressed in the parameter $\theta$, form
precisely the space
\[
 \mathcal T_t=\operatorname{span}_{\R}
 \{1,\cos\theta,\sin\theta,\ldots,\cos t\theta,\sin t\theta\}.
\]
Indeed, every coordinate of $\gamma(\theta)$ belongs to $\mathcal T_1$,
so each monomial of total degree at most $t$ restricts to an element of
$\mathcal T_t$. For the reverse inclusion, the affine linear functions
\[
 a(x)=\frac{(x-c,u)}r,\qquad b(x)=\frac{(x-c,v)}r
\]
restrict to $\cos\theta$ and $\sin\theta$, respectively. The real and
imaginary parts of $(a(x)+ib(x))^k$ restrict to $\cos k\theta$ and
$\sin k\theta$ for $0\leq k\leq t$.

A nonzero element of $\mathcal T_t$ has at most $2t$ distinct zeros
modulo $2\pi$: after writing it as a Laurent polynomial in
$z=e^{i\theta}$ and multiplying by $z^t$, one obtains a nonzero
polynomial of degree at most $2t$. Consequently, evaluation at any
$2t+1$ distinct points is an isomorphism from $\mathcal T_t$ to
$\R^{2t+1}$. The $2t+2$ evaluation functionals therefore have
a one-dimensional space of real linear relations, and every nonzero
relation has all its coefficients nonzero.

To determine the coefficients, put $N=2t+2$, $z_j=e^{i\theta_j}$, and
$P(z)=\prod_{j=1}^N(z-z_j)$. Comparing the coefficients of $z^{N-1}$
in the Lagrange interpolation formula for $z^q$ gives
\begin{equation}
\label{eq:minimal-null-lagrange-identity}
 \sum_{j=1}^N\frac{z_j^q}{P'(z_j)}=
 \begin{cases}
  0,&0\leq q\leq N-2,\\
  1,&q=N-1.
 \end{cases}
\end{equation}
In particular, $(z_j^t/P'(z_j))_{j=1}^N$ annihilates all the Laurent
monomials $z^k$, $-t\leq k\leq t$. If $S=\sum_{j=1}^N\theta_j$, then
\[
 P'(z_j)=(2i)^{N-1}e^{iS/2}z_j^t
          \prod_{k\ne j}\sin\frac{\theta_j-\theta_k}{2}.
\]
Multiplication by the common factor $(2i)^{N-1}e^{iS/2}$ gives
real, nonzero coefficients. Since the relation space is one-dimensional,
all null weight functions are given by
\begin{equation}
\label{eq:minimal-spherical-null-weights}
 w(x_j)=\beta\left(\prod_{k\ne j}
              \sin\frac{\theta_j-\theta_k}{2}\right)^{-1},
 \qquad \beta\in\R_{\ne0}.
\end{equation}

The chord-length identity gives
\[
 \|x_j-x_\ell\|=2r\left|\sin\frac{\theta_j-\theta_\ell}{2}\right|.
\]
Put $D_j=\prod_{\ell\ne j}\|x_j-x_\ell\|$.
There are $N-j$ negative factors in the sine product for $x_j$.
Since $N$ is even, we obtain
\[
 \prod_{\ell\ne j}\sin\frac{\theta_j-\theta_\ell}{2}
       =(-1)^j\frac{D_j}{(2r)^{N-1}}.
\]
Thus \eqref{eq:minimal-spherical-null-weights} is equivalent to
\eqref{eq:minimal-null-distance-weights} with $\alpha=-(2r)^{2t+1}\beta$.
Every $D_j$ is strictly positive, so
\eqref{eq:minimal-null-distance-weights} implies cyclic sign alternation
and $|X^+|=|X^-|=t+1$.

Finally, for the weights in \eqref{eq:minimal-spherical-null-weights},
identity \eqref{eq:minimal-null-lagrange-identity} with $q=N-1$ yields
\[
 \sum_{j=1}^N w(x_j)z_j^{t+1}
       =\beta(2i)^{N-1}e^{iS/2}\ne0.
\]
At least one of the real and imaginary parts of
$(a(x)+ib(x))^{t+1}$ therefore has a nonzero weighted sum.
Hence the null design has strength exactly $t$.
\end{proof}

The alternation of signs can also be obtained from the discrete
zero theorem for convex polygonal lines
\cite[Theorem~6]{Musin1998}, applied in the plane of the circle.

\begin{cor}
\label{cor:minimal-null-equal-weights}
Let $(X,w)$ be a spherical null $t$-design on $\bS^{n-1}$,
where $n\ge2$, $t\ge1$, and $|X|=2t+2$.
Suppose that $|w(x)|=a$ for every $x\in X$, where $a>0$.
Then $X^+$ and $X^-$ are two disjoint regular $(t+1)$-gons on the
unique circle containing $X$. For $t=1$, a regular $2$-gon means
an antipodal pair with respect to the centre of this circle.

Conversely, any two disjoint regular $(t+1)$-gons on the same
nondegenerate circle in $\bS^{n-1}$, assigned weights $a$ and
$-a$, form a spherical null design of strength exactly $t$ and
minimum size $2t+2$.
\end{cor}

\begin{proof}
Put $m=t+1$. By
Theorem~\ref{thm:minimal-spherical-null-classification}, the two sign
classes each have $m$ points and lie on a common circle.
Use angular coordinates to identify this circle with the unit circle in $\mathbb C$,
and let $z_1,\ldots,z_m$ and $u_1,\ldots,u_m$ represent $X^+$ and
$X^-$, respectively. The null conditions give
\[
 \sum_{j=1}^m z_j^k=\sum_{j=1}^m u_j^k
 \qquad(1\le k\le m-1).
\]
Write
\[
 \begin{aligned}
 F(z)&=\prod_{j=1}^m(z-z_j)=\sum_{k=0}^m c_kz^{m-k},\\
 G(z)&=\prod_{j=1}^m(z-u_j)=\sum_{k=0}^m d_kz^{m-k},
 \qquad c_0=d_0=1.
 \end{aligned}
\]
By Newton's identities, $c_k=d_k$ for $1\le k\le m-1$.
Moreover, $c_m\ne d_m$: otherwise $F=G$, contradicting the
fact that the two nonempty root sets are disjoint.

At this point, $F/(c_m-d_m)$ and $G/(c_m-d_m)$ satisfy the
hypotheses of the polynomial characterization in
\cite[Section~3.1, Problem~7]{BhattacharyaChenTSTST2020}.
We give the coefficient argument below for completeness.

All roots of $F$ have modulus one. Thus $|c_m|=1$, and comparison
of the roots and leading coefficients gives the polynomial identity
\[
 F(z)=c_m z^m\overline{F(1/\overline z)}.
\]
Here the right-hand side is interpreted by expanding it as a
polynomial, so the identity also holds at $z=0$.
Comparing coefficients yields
\[
 c_k=c_m\overline{c_{m-k}}\qquad(0\le k\le m).
\]
The same argument gives $d_k=d_m\overline{d_{m-k}}$.
For $1\le k\le m-1$, using the common intermediate coefficients,
we obtain
\[
 0=(c_m-d_m)\overline{c_{m-k}}.
\]
Since $c_m\ne d_m$, all intermediate coefficients vanish.
It follows that $F(z)=z^m+c_m$ and $G(z)=z^m+d_m$.
Thus their roots form regular $m$-gons.

Conversely, write the normalized root sets of two disjoint regular
$m$-gons as $\{\zeta\omega^j:0\le j<m\}$ and
$\{\eta\omega^j:0\le j<m\}$, where $|\zeta|=|\eta|=1$ and
$\omega=e^{2\pi i/m}$. For $1\le k\le m-1$, both power sums
vanish by the geometric-series identity. The constant sums both
equal $m$, and the negative power sums are complex conjugates of
the positive ones. The signed difference therefore annihilates
$\mathcal T_{m-1}$. By the restriction argument in the preceding
theorem, it annihilates every ambient polynomial of degree at
most $t$. Exact strength follows from that theorem, and minimum
size follows from the spherical lower bound $|X|\ge2t+2$.
\end{proof}

Put $m=t+1$. After choosing the angular origin and interchanging the two sign
classes if necessary, the angles can be written as
\[
 X^+:\quad \frac{2\pi j}{m},\qquad
 X^-:\quad \delta+\frac{2\pi j}{m}
 \quad(0\le j<m),\qquad 0<\delta<\frac{2\pi}{m}.
\]
Thus the relative rotation is arbitrary subject to disjointness.
The union need not be a regular $2m$-gon: this occurs precisely
when $\delta=\pi/m$.

Under the assumptions of
Corollary~\ref{cor:polynomial-representation} in the spherical case,
let $t\ge1$ and let $X\subset\Omega$ have $2t+2$ points.
By Proposition~\ref{prop:design-polynomial-space}, $X$ supports a
null $t$-design in $(\Omega,\rho)$ if and only if $\tau(X)$ lies on
a circle. The weights are obtained from
\eqref{eq:minimal-null-distance-weights} by replacing
$\|x_j-x_\ell\|$ with $\|\tau(x_j)-\tau(x_\ell)\|$.
The entire circle need not be contained in $\tau(\Omega)$.

\subsubsection{A regular-polygon example}

Bannai~\cite{Bannai2013} noted the existence of spherical null
$t$-designs with $2(t+1)$ points. An explicit example is obtained
by taking $n\ge2$, $t\ge1$, and
\[
 x_j=\left(\cos\frac{j\pi}{t+1},\sin\frac{j\pi}{t+1},
                 0,\ldots,0\right),\qquad
 w(x_j)=(-1)^j\quad(0\le j<2t+2).
\]
The even and odd vertices form two regular $(t+1)$-gons on a
common great circle, with relative rotation $\pi/(t+1)$.
Corollary~\ref{cor:minimal-null-equal-weights} gives a null design
of strength exactly $t$, attaining the bound in
Corollary~\ref{cor:Bannai-conjecture}. Here the union is a regular
$2(t+1)$-gon; for $t=1$, each sign class is an antipodal pair.

\subsection{Examples of minimum-size Euclidean null designs}

Consider $I=[-1,1]$ as a Euclidean polynomial space
with separation function $\rho(x,y)=xy$ and inner product
\[
    \langle f,g\rangle=\int_I fg\,d\mu,
    \qquad d\mu(x)=\frac12\,dx.
\]
For every integer $t\ge0$, choose distinct points
$a_1,\ldots,a_{t+2}\in I$, and put
\[
 X=\{a_1,\ldots,a_{t+2}\},\qquad
 w(a_j)=\frac{1}{\displaystyle\prod_{\ell\ne j}(a_j-a_\ell)}.
\]
These weights are all nonzero. For every polynomial $f$ of degree
at most $t$, the Lagrange interpolation formula gives
\[
 f(x)=\sum_{j=1}^{t+2} f(a_j)
       \frac{\displaystyle\prod_{\ell\ne j}(x-a_\ell)}
            {\displaystyle\prod_{\ell\ne j}(a_j-a_\ell)}.
\]
Comparing the coefficients of $x^{t+1}$ yields
\[
 \sum_{j=1}^{t+2}w(a_j)f(a_j)=0.
\]
Thus $(X,w)$ is a Euclidean null $t$-design with $t+2$ points.
Hence the universal lower bound $|X|\ge t+2$ for Euclidean
polynomial spaces is sharp for every $t\ge0$.

For odd $t$, a minimum-size null design also arises as the
difference of two positive quadrature formulas.
Let $m\ge1$ and $t=2m-1$, and let $(X_1,\lambda_1)$ and
$(X_2,\lambda_2)$ be the $m$-point and $(m+1)$-point
Gauss--Legendre formulas for $\mu$.
Their weights are positive and their node sets are disjoint,
since consecutive Legendre polynomials have strictly interlacing
zeros. On $X=X_1\sqcup X_2$, set $\lambda=\lambda_1$ on $X_1$
and $\lambda=-\lambda_2$ on $X_2$. Both formulas are exact
through degree $t$, so
\[
 \sum_{x\in X}\lambda(x)f(x)
 =\sum_{x\in X_1}\lambda_1(x)f(x)
  -\sum_{x\in X_2}\lambda_2(x)f(x)=0
 \qquad(f\in P_t(\R)).
\]
Thus $(X,\lambda)$ is a null $t$-design with the minimum number
$m+(m+1)=t+2$ of points.

\subsection{Existence of null designs}
In contrast to the explicit constructions above, dimension counting
ensures that every sufficiently large finite subset contains a
null support. It does not determine the support or
weights explicitly, nor their suitability for switching.

\begin{prop}
\label{prop:existence}
Let $(\Omega,\rho)$ be a polynomial space, and let
$X\subset\Omega$ be a finite subset.
If
\[
    |X|>\dim \Pol_\rho(\Omega,t),
\]
then there exist a nonempty subset $X'\subseteq X$ and a weight function
\[
    w:X'\to \R_{\ne0}
\]
such that
\[
    \sum_{x\in X'} w(x)f(x)=0
\]
for every $f\in\Pol_\rho(\Omega,t)$.
In particular, $(X',w)$ is a null $t$-design.
\end{prop}

\begin{proof}
For each $x\in X$, let
\[
    L_x:\Pol_\rho(\Omega,t)\to\mathbb{R},
    \qquad
    L_x(f)=f(x),
\]
be the evaluation functional at $x$.
Since
\[
    |X|>\dim\Pol_\rho(\Omega,t)
         =\dim\Pol_\rho(\Omega,t)^*,
\]
the family of evaluation functionals
\[
    \{L_x:x\in X\}\subset\Pol_\rho(\Omega,t)^*
\]
is linearly dependent.
Hence there exist real numbers $c_x$, not all zero, such that
\[
    \sum_{x\in X} c_xL_x=0.
\]
Set
\[
    X'=\{x\in X:c_x\neq 0\}
\]
and define $w:X'\to\R_{\ne0}$ by $w(x)=c_x$.
Then $X'$ is nonempty and, for every
$f\in\Pol_\rho(\Omega,t)$,
\[
    \sum_{x\in X'}w(x)f(x)
    =
    \sum_{x\in X}c_xL_x(f)
    =
    0.
\]
Thus $(X',w)$ is a null $t$-design.
\end{proof}

Equivalently, for a basis $f_1,\ldots,f_D$ of
$\Pol_\rho(\Omega,t)$, the evaluation matrix
$M=(f_j(x))_{x\in X,\,1\le j\le D}$ expresses the relation above
as $M^{\mathsf T}c=0$; the required support is
$X'=\{x\in X:c_x\ne0\}$.

For the standard sphere, the dimension bound in
Proposition~\ref{prop:existence} takes the following form.

\begin{cor}
\label{cor:spherical-null-existence}
Let $n\ge2$ and $t\ge0$, and let $X$ be a finite subset of
$\bS^{n-1}$. Put
\[
 h_i=\dim\Harm_i(\R^n),\qquad D_{n,t}=\sum_{i=0}^t h_i,
\]
where $\Harm_i(\R^n)$ denotes the space of homogeneous harmonic
polynomials of degree $i$.
If $|X|>D_{n,t}$, then there exist a nonempty subset $X'\subseteq X$
and a weight function $w:X'\to\R_{\ne0}$ such that $(X',w)$ is a
spherical null $t$-design.
\end{cor}

\begin{proof}
The restrictions to $\bS^{n-1}$ of polynomials of degree at most $t$
form the direct sum
\[
 \bigoplus_{i=0}^t\Harm_i(\R^n)\big|_{\bS^{n-1}},
\]
which has dimension $D_{n,t}$.
The result follows from
Proposition~\ref{prop:existence}.
\end{proof}

Proposition~\ref{prop:existence} guarantees a null design supported
on some nonempty subset of $X$ under its dimension hypothesis,
but does not ensure nonzero null weights at every point of $X$.
For combinatorial null designs, Deza et al.\ \cite[Theorem~6]{DezaFranklSinghi1983}
describe null supports as unions of circuits. This viewpoint leads
to the following criterion for a prescribed support in a general
polynomial space.

\begin{rem}
\label{rem:full-support-null}
Let $(\Omega,\rho)$ be a polynomial space, let $t\ge0$ be an
integer, and put $V=\Pol_\rho(\Omega,t)$.
For a nonempty finite set $X\subset\Omega$, define
$L_x\in V^*$ by $L_x(f)=f(x)$ and put $H_x=\ker L_x$.
Let $M_X$ be the matroid on $X$ represented by the indexed
family $(L_x)_{x\in X}$, with independence given by linear
independence of the functionals. An element $x$ is a coloop
when $L_x$ is not in the span of the other functionals.
Then the following conditions are equivalent:
\begin{enumerate}
\item[(i)] There is $w:X\to\R_{\ne0}$ such that
$\sum_{x\in X}w(x)f(x)=0$ for every $f\in V$.
\item[(ii)] The matroid $M_X$ has no coloops.
\item[(iii)] For every $x\in X$,
\[
 \bigcap_{y\in X\setminus\{x\}}H_y
 =\bigcap_{y\in X}H_y.
\]
\end{enumerate}
An intersection over the empty family means $V$.

Indeed, $V$ is finite-dimensional and contains $1$, so $L_x(1)=1$
and each $H_x$ is a hyperplane. Condition (i) expresses every
$L_x$ as a linear combination of the others, giving (ii).
Conversely, under (ii), the relation space
\[
 K=\left\{u\in\R^X:\sum_{x\in X}u(x)L_x=0\right\}
\]
contains, for each $x$, an element with nonzero $x$-coordinate.
Thus each coordinate kernel in $K$ is a proper hyperplane.
Avoiding their finite union over $\R$ yields a relation nonzero
in every coordinate, proving (i).
Finally, by finite-dimensional duality, two families of linear
functionals have the same common kernel if and only if they
span the same subspace of $V^*$. This gives (ii)$\Leftrightarrow$(iii).

The circuits of $M_X$ are exactly the inclusion-minimal null
supports contained in $X$, and (ii) says that $X$ is a union
of such circuits. The remark following Theorem~6 of
\cite{DezaFranklSinghi1983} also notes the validity of this
observation for general chain matroids.
These conditions do not prescribe the signs or magnitudes of
the null weights. The weight conditions for switching in
Section~\ref{sect:4} must be checked separately.
\end{rem}

\section{Switching of designs in spherical polynomial spaces}
\label{sect:4}

In this section, we study a switching operation for designs in
spherical polynomial spaces. To apply a null configuration to a given design, its placement
and weights must allow the modification without producing
negative weights.
Our starting point is the local modification argument used by
Sali~\cite[Lemma~2.3]{Sali1994} in his study of the rigidity of spherical
designs.
If a configuration on a small sphere becomes a spherical
$t$-design after centring and rescaling to the unit sphere,
rotating it about that centre within its affine span preserves
its ambient polynomial moments through degree $t$.
The difference between the original and rotated configurations
is therefore a null relation. The null condition requires only
that the two parts have matching weighted polynomial moments;
neither part need itself be a spherical $t$-design. We formulate
this replacement principle in spherical polynomial spaces.

We first describe the normalization of a small sphere, then
give complete and partial switching constructions.
Corollary~\ref{cor:circle-switching} applies the minimum-size
classification to switching on circles.

\begin{lem}
\label{lem:embed}
Let $(\Omega, \rho)$ be a spherical polynomial space of dimension $n$.
Suppose there exists an injection $\tau: \Omega \to \bS^{n-1}$ that is affinely compatible with $\rho$.
Let $Y\subset \Omega$, and suppose that there exist an
$m$-dimensional linear subspace $T^m\subset\mathbb R^n$ and
$z_0\in (T^m)^\perp$ such that
$\tau(Y)\subset (z_0+T^m)\cap \bS^{n-1}$.
Then
\[
\|\tau(y)-z_0\|^2=1-\|z_0\|^2 \qquad \text{for all } y\in Y.
\]
In particular, if $r=\sqrt{1-\|z_0\|^2}>0$,
then
\[
    \tau_Y(y):=\frac{\tau(y)-z_0}{r}
\]
defines an injective map $\tau_Y:Y\rightarrow \bS^{m-1}\subset T^m$.
\end{lem}

\begin{proof}
Since $\tau(y)\in z_0+T^m$, we have
$\tau(y)-z_0\in T^m$.  Since $z_0\in (T^m)^\perp$,
\[
    \bigl(z_0,\tau(y)-z_0\bigr)=0.
\]
Hence $1=\|\tau(y)\|^2 =\|z_0\|^2+\|\tau(y)-z_0\|^2$,
which gives
\[
    \|\tau(y)-z_0\|^2=1-\|z_0\|^2.
\]
Thus $\tau_Y(y)\in \bS^{m-1}\subset T^m$.
Since $\tau$ is injective, so is $\tau_Y$.
\end{proof}

For a nonempty finite set $Y$ as above and a weight function
$\lambda:Y\to\R_{\ne0}$, we say that
$(Y,\lambda)$ is a null $t$-design with respect to $\tau_Y$ if
\[
 \sum_{y\in Y}\lambda(y)g(\tau_Y(y))=0
 \qquad\text{for every }g\in P_t(\R^m),
\]
where $T^m$ is identified with $\R^m$ by orthonormal coordinates.

\begin{thm}[Switching theorem]
\label{thm:switching1}
Let $(\Omega, \rho)$ be a spherical polynomial space of dimension $n$.
Suppose there exists an injection $\tau: \Omega \to \bS^{n-1}$
that is affinely compatible with $\rho$.
Let $(X,w)$ be a weighted $t$-design in $(\Omega, \rho)$,
and let
\[
    Y^+\subset X,
    \qquad
    Y^-\subset\Omega\setminus X,
    \qquad
    Y=Y^+\sqcup Y^-.
\]
Suppose that $Y$ satisfies the assumptions of
Lemma~\ref{lem:embed} and that
$r=\sqrt{1-\|z_0\|^2}>0$.
Let $\tau_Y:Y\to \bS^{m-1}\subset T^m$
be the induced spherical embedding.
Suppose that there exist a weight function
$\lambda:Y\to\mathbb R$ and a constant $c>0$ such that
\[
    c \lambda(y) = w(y) \quad (y \in Y^+),
    \qquad
    \lambda(y)<0 \quad (y \in Y^-),
\]
and that $(Y,\lambda)$ is a null $t$-design with respect to $\tau_Y$.
Then $(X',w')$ is a weighted $t$-design in $(\Omega,\rho)$,
where
\[
    X'=(X\setminus Y^+)\sqcup Y^-
\]
and $w':X'\to\mathbb R_{>0}$ is defined by
\[
    w'(x)=
    \begin{cases}
        w(x), & x\in X\setminus Y^+,\\
        -c \lambda(x), & x\in Y^-.
    \end{cases}
\]
\end{thm}

\begin{proof}
By Proposition~\ref{prop:design-polynomial-space}, it suffices
to verify the corresponding polynomial identities.

Since $(Y,\lambda)$ is a null $t$-design with respect to $\tau_Y$,
the definitions of $\lambda$ and $w'$ imply that,
for any $g \in P_t(\R^m)$,
\begin{align*}
\sum_{y \in Y^+} w(y) \; g \circ \tau_Y(y)
&= \sum_{y \in Y^+} c \lambda(y) \; g \circ \tau_Y(y)
\\
&= \sum_{y \in Y^-} (-c \lambda(y)) \; g \circ \tau_Y(y) \\
&= \sum_{y \in Y^-} w'(y) \; g \circ \tau_Y(y).
\end{align*}

Choose orthonormal coordinates on $\R^n$ such that
\[
T^m=\{0\}^{n-m}\times\R^m,
\qquad z_0=(a_1,\ldots,a_{n-m},0,\ldots,0).
\]
For any $f\in P_t(\R^n)$, there exists
$\widetilde g\in P_t(\R^m)$ such that
\[
f(a_1,\ldots,a_{n-m}, r x_{n-m+1}, \ldots, r x_n)
= \widetilde{g}(x_{n-m+1},\ldots,x_n).
\]
Hence,
\begin{align*}
\sum_{y \in Y^+} w(y) \; f \circ \tau(y)
&= \sum_{y \in Y^+} w(y) \; f(r \tau_Y(y) + z_0) \\
&= \sum_{y \in Y^+} w(y) \; \widetilde{g} \circ \tau_Y(y) \\
&= \sum_{y \in Y^-} w'(y) \; \widetilde{g} \circ \tau_Y(y) \\
&= \sum_{y \in Y^-} w'(y) \; f(r \tau_Y(y) + z_0) \\
&= \sum_{y \in Y^-} w'(y) \; f \circ \tau(y).
\end{align*}
Since $(X,w)$ is a weighted $t$-design, it follows that
\begin{align*}
\langle 1, f \circ \tau \rangle
&= \sum_{x \in X \setminus Y^+} w(x) \; f \circ \tau(x)
+ \sum_{y \in Y^+} w(y) \; f \circ \tau(y) \\
&= \sum_{x \in X \setminus Y^+} w(x) \; f \circ \tau(x)
+ \sum_{y \in Y^-} w'(y) \; f \circ \tau(y) \\
&= \sum_{y \in (X \setminus Y^+)\sqcup Y^-}
w'(y) \; f \circ \tau(y)
\qquad\text{for all }f \in P_t(\R^n).
\end{align*}
By Proposition~\ref{prop:design-polynomial-space},
$(X',w')$ is a weighted $t$-design in $(\Omega,\rho)$.
\end{proof}

Theorem~\ref{thm:switching1} removes the entire weight at each
point of $Y^+$. The same argument also applies when only part of
the weight at each point of $Y^+$ is removed, so that the points
in $Y^+$ remain in the support.
We refer to this operation as partial switching.
\begin{prop}
\label{prop:switching2}
Let $(\Omega, \rho)$ be a spherical polynomial space of dimension $n$.
Suppose there exists an injection $\tau: \Omega \to \bS^{n-1}$ that is affinely compatible with $\rho$.
Let $(X,w)$ be a weighted $t$-design in $(\Omega, \rho)$,
and let
\[
    Y^+\subset X,
    \qquad
    Y^-\subset\Omega\setminus X,
    \qquad
    Y=Y^+\sqcup Y^-.
\]
Suppose that $Y$ satisfies the assumptions of Lemma~\ref{lem:embed}
and that $r=\sqrt{1-\|z_0\|^2}>0$.
Let $\tau_Y:Y\to \bS^{m-1}\subset T^m$
be the induced spherical embedding.
Suppose that there exists a weight function
$\lambda:Y\to\mathbb R$ such that
\[
    0 < \lambda (y) < w(y), \quad y \in Y^+,
    \qquad
    \lambda(y) <0, \quad y \in Y^-,
\]
and $(Y,\lambda)$ is a null $t$-design with respect to $\tau_Y$.
Then $(X'',w'')$ is a weighted $t$-design in $(\Omega,\rho)$,
where $X''=X\sqcup Y^-$ and
$w'':X''\to\R_{>0}$ is defined by
\[
    w''(x)=
    \begin{cases}
        w(x), & x\in X\setminus Y^+,\\
        w(x) - \lambda (x), & x \in Y^+, \\
        - \lambda(x), & x\in Y^-.
    \end{cases}
\]
\end{prop}

\begin{proof}
As in the proof of Theorem~\ref{thm:switching1}, the affine map
$u\mapsto z_0+ru$ preserves the degree bound, so
$\sum_{y\in Y}\lambda(y)f(\tau(y))=0$ for every $f\in P_t(\R^n)$.
Hence
\[
 \sum_{x\in X''}w''(x)f(\tau(x))
 =\sum_{x\in X}w(x)f(\tau(x))
   -\sum_{y\in Y}\lambda(y)f(\tau(y))
 =\langle1,f\circ\tau\rangle.
\]
By the inequalities on $\lambda$, all weights of $w''$ are
positive. The result follows from
Proposition~\ref{prop:design-polynomial-space}.
\end{proof}

The alternating signs and explicit weights in the minimum-size
classification now give a concrete partial switching on a circle.
The selected points of the original design may have arbitrary
positive weights.

\begin{cor}
\label{cor:circle-switching}
Let $n\ge2$ and $t\ge1$, and let $(Z,w)$ be a positive weighted
spherical $t$-design on $\bS^{n-1}$.
Suppose that $A\subset Z$ consists of $t+1$ distinct points on
a nondegenerate circle $C\subset\bS^{n-1}$.
Choose one point of $C\setminus Z$ in each open arc between
consecutive points of $A$, and let $B$ be the set of these points.
List $A\sqcup B=\{x_1,\ldots,x_{2t+2}\}$ in cyclic order,
starting with a point of $A$, and put
\[
 \lambda(x_j)=
 \frac{(-1)^{j-1}}
 {\displaystyle\prod_{\ell\ne j}\|x_j-x_\ell\|}.
\]
Then $\lambda$ is positive on $A$ and negative on $B$.
Set
\[
 c_* = \min_{a\in A}\frac{w(a)}{\lambda(a)}.
\]
For every $0<c<c_*$, the weight function
\[
 w_c(x)=
 \begin{cases}
  w(x), & x\in Z\setminus A,\\
  w(x)-c\lambda(x), & x\in A,\\
  -c\lambda(x), & x\in B
 \end{cases}
\]
makes $(Z\sqcup B,w_c)$ a positive weighted spherical $t$-design.
\end{cor}

\begin{proof}
Since $Z$ is finite, we may choose a point of $C\setminus Z$
in each of the indicated open arcs. The sets $A$ and $B$
alternate in cyclic order. By
Theorem~\ref{thm:minimal-spherical-null-classification},
$(A\sqcup B,\lambda)$ is a spherical null $t$-design with
the stated signs. The choice of $c$ makes every value of $w_c$
strictly positive. Moreover, for every $f\in P_t(\R^n)$,
\[
 \sum_{x\in Z\sqcup B}w_c(x)f(x)
 =\sum_{x\in Z}w(x)f(x)
   -c\sum_{x\in A\sqcup B}\lambda(x)f(x)
 =\sum_{x\in Z}w(x)f(x).
\]
Thus all moments through degree $t$ are preserved.
\end{proof}

The same formula removes all points of $A$ at once precisely
when $w(a)/\lambda(a)$ is independent of $a\in A$.
Taking $c$ equal to this common value and deleting the zero
weights gives a positive weighted $t$-design on
$(Z\setminus A)\sqcup B$.
This is the complete replacement described in
Theorem~\ref{thm:switching1}, whereas $0<c<c_*$ gives
partial switching.

For a spherical polynomial space with a compatible embedding
$\tau$, the same construction applies when the selected points
of $A\sqcup B$ lie in $\tau(\Omega)$; the existence of a point
of $\tau(\Omega)$ in each open arc is an additional requirement.

\section{Trade lifting and examples of design switching}
\label{sect:5}

We begin with the cube example, then establish the spherical
realizations of Hamming and Johnson trades and apply them to
the binary OA and Pasch examples. We also explain their common
parity construction and conclude with an equal-weight switching
between two noncongruent $240$-point spherical $5$-designs.

Throughout this section, Hamming and Johnson spaces are equipped
with the uniform probability measure, while spheres are equipped
with normalized surface measure.
Integer multiplicities are denoted by $m$ and are normalized to weights
$w(x)=m(x)/M$, where $M=\sum_x m(x)$. We write $\mathbf{1}_E$ for
the indicator of a set $E$ and $\mathbf{1}$ for the all-ones vector.

\subsection{Geometric switching in the sphere}

\begin{exam}
\label{exam:cube}
Let $X$ be the eight vertices of a regular cube inscribed in $\bS^2$,
with $w(x)=1/8$ for every $x\in X$. Then $(X,w)$ is a spherical
$3$-design. Let $Y^+\subset X$ be the four vertices of one face,
and rotate $Y^+$ about the centre of that face in its affine plane
to obtain $Y^-$. Choose a rotation angle not congruent to $0$
modulo $\pi/2$, so that $Y^-\cap X=\varnothing$.

On the circle determined by the face, $Y^+$ and $Y^-$ are disjoint
regular squares. By Corollary~\ref{cor:minimal-null-equal-weights},
their signed difference, with $\lambda=1$ on $Y^+$ and
$\lambda=-1$ on $Y^-$, is a null $3$-design, also after
normalizing the circle to the unit circle.
The null support $Y^+\sqcup Y^-$ has affine dimension $2$
and the minimum size $8=2(3+1)$, although the cube itself
has affine dimension $3$.
Theorem~\ref{thm:switching1}, with $c=1/8$, therefore gives a
spherical $3$-design $(X',w')$, where
\[
 X'=(X\setminus Y^+)\sqcup Y^-,
 \qquad w'(x)=\frac18\quad(x\in X').
\]
Thus four vertices are replaced while all remaining vertices
and the equal weights are retained.
\end{exam}

\subsection{Trade lifting for the standard embeddings}

We use the notation of Section~\ref{sect:2}.
Let $\mathcal A$ be a set of size $q\ge2$, let $d\ge1$, and write
$[m]=\{1,\ldots,m\}$ for every positive integer $m$. Equip $H(d,q)=\mathcal A^d$ with Hamming distance
$\rho_H$. Its polynomial spaces have the description
\[
 \Pol_{\rho_H}(H(d,q),r)=
 \Span_{\R}\left\{
   \prod_{i\in I}\mathbf{1}_{\{x_i=a_i\}}:
   I\subseteq[d],\ |I|\leq r,\ a_i\in\mathcal A
 \right\}.
\]
Choose the vertices $\{e_a:a\in\mathcal A\}$ of a
regular simplex in $\R^{q-1}$, normalized by
\[
 (e_a,e_b)=
 \begin{cases}
   1,&a=b,\\
   -1/(q-1),&a\ne b.
 \end{cases}
\]
The standard spherical embedding is
\[
 \tau_H(x)=\frac{1}{\sqrt d}
       (e_{x_1},\ldots,e_{x_d})
 \in\bS^{d(q-1)-1}.
\]

The embedding satisfies
\[
 (\tau_H(x),\tau_H(y))
 =1-\frac{q}{d(q-1)}\rho_H(x,y).
\]
For $1\le k\le v-1$, equip the Johnson space
$J(v,k)=\binom{[v]}{k}$ with
$\rho_J(A,B)=k-|A\cap B|$. Its polynomial spaces have the description
\[
 \Pol_{\rho_J}(J(v,k),r)=
 \Span_{\R}\left\{
   \mathbf{1}_{\{D\subseteq A\}}:
   D\subseteq[v],\ |D|\leq r
 \right\}.
\]
Put $L=\mathbf{1}^{\perp}\subset\R^v$ and identify $L$ with
$\R^{v-1}$ by an orthonormal basis. Its unit sphere is thereby
identified with $\bS^{v-2}$. The standard spherical embedding is
\[
 \tau_J(A)=
 \frac{\mathbf{1}_A-\frac{k}{v}\mathbf{1}}
      {\sqrt{k(v-k)/v}}
 \in L,\qquad \|\tau_J(A)\|=1.
\]
In particular,
\begin{equation}
 \label{eq:johnson-inner-product}
 (\tau_J(A),\tau_J(B))
 =\frac{v|A\cap B|-k^2}{k(v-k)}
 =1-\frac{v}{k(v-k)}\rho_J(A,B).
\end{equation}

Thus both embeddings are affinely compatible with their respective
separation functions. The polynomial-space dimensions are $d(q-1)$
and $v-1$, respectively.

For either space $\Omega$, represent a trade by nonnegative
integer-valued multiplicity functions $m^+,m^-$ on $\Omega$.
After cancelling common occurrences, assume that their supports
$T^+,T^-$ are disjoint and nonempty.
The sets $T^+$ and $T^-$ are called the positive and negative
legs of the trade, respectively.
Put
\[
 T=T^+\sqcup T^-,\qquad
 \lambda(x)=m^+(x)-m^-(x)\quad(x\in T).
\]
Thus $\lambda:T\to\R_{\ne0}$ is positive on $T^+$ and negative on
$T^-$. We call the trade a $t$-trade if $(T,\lambda)$ is a null
$t$-design in $(\Omega,\rho)$. For a simple trade,
$m^\pm=\mathbf{1}_{T^\pm}$.
The common value $\sum_{x\in T^+}m^+(x)=\sum_{x\in T^-}m^-(x)$
is the \textit{volume} of the trade; it counts multiplicities on
one leg rather than distinct points in the full support.

\begin{prop}
\label{prop:trade-lifting}
Let
\[
    (\Omega,\rho)=(H(d,q),\rho_H)
    \quad
    (\text{resp. }(\Omega,\rho)=(J(v,k),\rho_J)),
\]
and let $\tau=\tau_H$ (resp.\ $\tau=\tau_J$) be the standard
embedding into $\bS^{N-1}$, where
\[
    N=d(q-1)
    \quad
    (\text{resp. }N=v-1).
\]
For every $r\geq0$,
\[
    \Pol_\rho(\Omega,r)
    =
    \{f\circ\tau:\Omega\to\R
      \mid f\in\R_r[x_1,\ldots,x_N]\}.
\]
If $(T,\lambda)$ is a trade of strength $t$, then
$(\tau(T),\widehat\lambda)$ is a spherical null $t$-design, where
\[
 \widehat\lambda:\tau(T)\to\R_{\ne0},
 \qquad \widehat\lambda(\tau(x))=\lambda(x).
\]
If the trade has exact strength $t$, its spherical lift also has
exact strength $t$.
\end{prop}

\begin{proof}
The polynomial representation follows from
Corollary~\ref{cor:polynomial-representation}. 
In the Hamming embedding, each ambient coordinate depends on a single entry of the word, while
\[
 \mathbf{1}_{\{x_i=a\}}
 =\frac{1+(q-1)(e_{x_i},e_a)}{q}
\]
is affine linear in the $i$th simplex block. Products of at most
$r$ indicators therefore span precisely the restrictions of
polynomials of degree at most $r$.

For the Johnson embedding, put $z_i(A)=\mathbf{1}_{\{i\in A\}}$.  Since
\[
 z_i(A)=\frac{k}{v}+\sqrt{\frac{k(v-k)}{v}}\,
                  \tau_J(A)_i,
\]
the incidence coordinates and the spherical coordinates are related
affinely.  After using $z_i^2=z_i$ on $\{0,1\}^v$, every polynomial of
degree at most $r$ is a linear combination of the square-free
monomials
$\prod_{i\in D}z_i=\mathbf{1}_{\{D\subseteq A\}}$ with $|D|\leq r$.
Polynomials on $L$ extend to $\R^v$ without increasing their degree,
so this also applies to the intrinsic coordinates on $L$.
This proves the asserted descriptions of the polynomial spaces.
The assertions on null designs and exact strength follow.
\end{proof}

When $r$ equals the diameter ($d$ and $\min\{k,v-k\}$ for the Hamming and Johnson cases,
respectively), $\Pol_\rho(\Omega,r)$
contains all functions on $\Omega$.

\begin{cor}
\label{cor:local-trade-switching}
Let $(Z,w)$ be a weighted spherical $t$-design on $\bS^{n-1}$,
and let $(\tau(T),\widehat\lambda)$ be a lifted trade as above.
Put $V=\operatorname{span}\tau(T)$.
Suppose that a linear isometry $U:V\to\R^n$ satisfies
$U\tau(T^+)\subseteq Z$.
Extend $w$ by zero outside $Z$, writing $\bar w$ for this extension,
and define
\[
 \Lambda(z)=
 \begin{cases}
  \lambda(x),&z=U\tau(x),\ x\in T,\\
  0,&z\notin U\tau(T).
 \end{cases}
\]
If $c>0$ satisfies $\bar w(z)-c\Lambda(z)\ge0$ for every
$z\in\bS^{n-1}$, then $(Z',w')$ is a weighted spherical
$t$-design, where
\[
 Z'=\{z\in Z\cup U\tau(T):\bar w(z)-c\Lambda(z)>0\},
 \qquad
 w'(z)=\bar w(z)-c\Lambda(z)\quad(z\in Z').
\]
\end{cor}

\begin{proof}
For $f\in\R_t[x_1,\ldots,x_n]$, the function $f\circ U$ on $V$
extends to an ambient polynomial of degree at most $t$.
Thus
\[
 \sum_{x\in T}\lambda(x)f(U\tau(x))=0.
\]
Subtracting $c$ times this identity from the design identity for
$(Z,w)$ proves the assertion. Discarding zero weights ensures that
$w':Z'\to\R_{>0}$.
\end{proof}

Both legs must be placed using the same isometry $U$.
Matching the Gram matrix of the positive leg alone does not justify
an arbitrary placement of the negative leg.
This operation allows some weight on the positive leg to remain,
and the negative leg may meet $Z$.

\subsection{The binary OA trade}

We write $OA(N,d,q,t)$ for an $N$-row, $d$-column array over
an alphabet $\mathcal A$ of size $q$ in which every $t$-tuple occurs equally often in
every choice of $t$ columns; repeated rows are allowed.
See \cite{HedayatSloaneStufken1999}.
An OA trade of strength $t$ is a pair
of multisets in $\mathcal A^d$ having the
same projections, with multiplicity, onto every set of at most $t$
coordinates.  Equivalently,
the multiplicity difference $\lambda=m^+-m^-$ on its support
annihilates $\Pol_{\rho_H}(H(d,q),t)$.  For the relation between
binary orthogonal arrays and set-theoretic trades, see~\cite{GhorbaniKamaliKhosrovshahiKrotov2020}.

\begin{exam}
\label{exam:binary-oa}
In $H(3,2)$ consider
\[
 T^+=\{000,011,101,110\},\qquad
 T^-=\{001,010,100,111\}.
\]
Both legs are $OA(4,3,2,2)$, so $(T^+,T^-)$ is an OA trade of
strength $2$.  Use the binary spherical embedding
\[
 \tau_H(x_1x_2x_3)
   =\frac{1}{\sqrt3}
      \bigl((-1)^{x_1},(-1)^{x_2},(-1)^{x_3}\bigr).
\]
Since distinct words in $T^+$ have distance $2$, their spherical
inner product is $-1/3$. Thus the positive leg is a regular
tetrahedron with Gram matrix
\begin{equation}
 \label{eq:tetrahedron-gram}
 G=\frac43I_4-\frac13J_4,
\end{equation}
where $J_4$ is the all-ones matrix.
Bitwise complementation interchanges the legs and acts antipodally,
so $\tau_H(T^-)=-\tau_H(T^+)$.

Let $T=T^+\sqcup T^-$, with $\lambda=1$ on $T^+$ and
$\lambda=-1$ on $T^-$. Proposition~\ref{prop:trade-lifting} shows
that $(\tau_H(T),\widehat\lambda)$ is a spherical null $2$-design.  Its strength is
exactly $2$: for
$f(y)=y_1y_2y_3$,
\[
 \sum_{x\in T^+}f(\tau_H(x))
 -\sum_{x\in T^-}f(\tau_H(x))
 =\frac{8}{3\sqrt3}\ne0.
\]
Consequently, switching replaces a regular tetrahedron by its
antipodal tetrahedron while preserving all moments through degree
$2$.

For partial switching, start with two copies of the positive leg:
\[
 m=2\mathbf{1}_{T^+},\qquad
 m'=m-\mathbf{1}_{T^+}+\mathbf{1}_{T^-}
    =\mathbf{1}_{\{0,1\}^3}.
\]
Both multiplicities total $8$. After normalization, weights $1/4$
on a tetrahedron become weights $1/8$ on the cube.
The initial spherical $2$-design thus becomes a spherical $3$-design:
switching preserves the prescribed strength and may increase it.
\end{exam}

\subsection{The Pasch trade in the Johnson scheme}

For a Johnson $t$-trade with $0\le t\le k$, the two legs contain each $t$-subset equally
often.  Equality for all subsets of size at most $t$ then follows, so
the signed difference annihilates $\Pol_{\rho_J}(J(v,k),t)$.  We use the following
standard volume-four trade.  For the classical theory of
$(v,k,t)$-trades and their use in switching designs, see
\cite{Hwang1986}; see also \cite{GhorbaniKamaliKhosrovshahiKrotov2020} for a recent
treatment of small trades.

\begin{exam}
\label{exam:pasch}
Let the ground set be $V=\{0,1,2,3,4,5\}$, and abbreviate
$\{0,2,3\}$ to $023$.  Put
\[
 T^+=\{023,045,124,135\},\qquad
 T^-=\{024,035,123,145\}.
\]
Each pair of points occurs equally often in the two legs.
Thus we obtain a $2$-trade in $J(6,3)$.

Embed $J(6,3)$ in the unit sphere of
$\mathbf{1}^{\perp}\subset\R^6$, identified with $\bS^4$, by
\[
 \tau_J(A)=\sqrt{\frac23}
       \left(\mathbf{1}_A-\frac12\mathbf{1}\right).
\]
Distinct blocks in $T^+$ meet in one point, so
\eqref{eq:johnson-inner-product} gives the tetrahedral Gram matrix
\eqref{eq:tetrahedron-gram}. Although the embedding lies in $\bS^4$,
the positive leg spans a three-dimensional subspace.
Complementation takes $T^+$ to $T^-$ and, since $k=v/2$,
$\tau_J(A^c)=-\tau_J(A)$. Thus the two legs have the same
antipodal tetrahedral geometry as Example~\ref{exam:binary-oa}.
Proposition~\ref{prop:trade-lifting} gives a spherical null
$2$-design with $\lambda=1$ on $T^+$ and $\lambda=-1$ on $T^-$.
This null design has strength exactly $2$: the cubic incidence monomial
$\mathbf{1}_{\{023\subseteq A\}}$ has signed sum $1$ against
$\lambda$ and is the restriction of an ambient polynomial of
degree at most $3$.

For partial switching, start with multiplicity $m(B)=2$ on every
block of $\binom V3$, a $2$-$(6,3,8)$ design since each pair
belongs to four blocks. The update
$m'=m-\mathbf{1}_{T^+}+\mathbf{1}_{T^-}$ has multiplicities
$1$ on $T^+$, $3$ on $T^-$, and $2$ elsewhere, and is again a
$2$-$(6,3,8)$ design. Both total multiplicities are $40$;
the normalized weights are $w(B)=2/40$ and $w'(B)=m'(B)/40$.

The embedded configuration $\tau_J(J(6,3))$, equipped with
uniform weights, is a spherical $2$-design on the unit sphere
in $L=\mathbf{1}^{\perp}$.
Indeed, if $P_L$ is the orthogonal
projection onto $L$, then
\[
 \frac1{20}\sum_{A\in\binom V3}\tau_J(A)=0,\qquad
 \frac1{20}\sum_{A\in\binom V3}\tau_J(A)\tau_J(A)^{\mathsf T}
 =\frac15P_L.
\]
These agree with the first and second moments of the sphere.
Hence, by Corollary~\ref{cor:local-trade-switching}, the normalized
switched weights also give a spherical
$2$-design on this sphere.
The tetrahedral legs themselves are spherical $2$-designs on their
own two-dimensional sphere, not on the ambient $\bS^4$;
their signed difference nevertheless satisfies the ambient null
$2$-design condition.
\end{exam}

\begin{rem}
\label{rem:parity-family}
Examples~\ref{exam:binary-oa} and~\ref{exam:pasch} are the $t=2$
case of a common construction. For $t\ge1$, put $d=t+1$ and
divide $\{0,1\}^d$ according to the parity of $\sum_i x_i$.
The two classes form an OA $t$-trade, since fixing at most $t$
coordinates leaves equally many extensions of each parity.
Their difference has strength exactly $t$, as the product
$\prod_{i=1}^d(-1)^{x_i}$ has nonzero sum against the parity weights.

For the Johnson realization, partition a $2d$-set into pairs
$\{a_i,b_i\}$ and associate $x$ with the transversal
$A_x=\{a_i:x_i=0\}\cup\{b_i:x_i=1\}$.
The parity classes give the standard minimum-volume $t$-trade
in $J(2d,d)$; see \cite[Lemma~4]{GhorbaniKamaliKhosrovshahiKrotov2020}.
With $f_i=(\mathbf e_{a_i}-\mathbf e_{b_i})/\sqrt2$, where
$\mathbf e_a$ denotes a standard coordinate vector, we have
\[
 \tau_J(A_x)=\frac1{\sqrt d}\sum_{i=1}^d(-1)^{x_i}f_i.
\]
Since the $f_i$ are orthonormal, this is the binary Hamming
embedding in these coordinates. For $t=2$, the pairs
$\{0,1\},\{2,5\},\{3,4\}$ give Example~\ref{exam:pasch}.

The common spherical image is a null $t$-design by
Proposition~\ref{prop:trade-lifting}. Its individual legs need
not be spherical $t$-designs even on the sphere in their linear
span: for $t\ge4$, their coordinate fourth moment is $1/d^2$,
whereas the spherical moment is $3/(d(d+2))$.
\end{rem}

\subsection{Complete switching in a spherical $5$-design}
\label{subsect:240point-switching}

We next combine a classical combinatorial $2$-trade with Victoir's
construction to illustrate Theorem~\ref{thm:switching1}, using
sign averaging to obtain a spherical null $5$-design.
The resulting replacement preserves equal weights and the
number of points, and produces a noncongruent spherical
$5$-design. Neither exchanged part is itself a spherical $5$-design
on the unit sphere of its linear span.

\begin{exam}
\label{exam:240point-switching}
Let $\mathcal B$ consist of the edge sets of
the $4$-cycles of $K_5$, with the edges labelled by
\[
 (1,\ldots,10)=(ab,ac,ad,ae,bc,bd,be,cd,ce,de).
\]
Each edge lies in six blocks, and each pair of edges lies in two.
Thus $\mathcal B$ is a $2$-$(10,4,2)$ design with $15$ blocks.
For $B\in\binom{[10]}4$, put
\[
 Q_B=\left\{\frac12\sum_{i\in B}\varepsilon_i e_i:
                  \varepsilon_i\in\{-1,1\}\right\}\subset\bS^9,
 \qquad X=\bigcup_{B\in\mathcal B}Q_B,
\]
where $e_1,\ldots,e_{10}$ is the standard orthonormal basis.

The full signed-permutation orbit of
$(1/2,1/2,1/2,1/2,0,\ldots,0)$ has
$2^4\binom{10}{4}=3360$ points and is a spherical $5$-design,
which arises from a particular case of the classical K\"ursch\'ak identity discussed in
\cite[Section~6]{NozakiSawa2014}.
Selecting the supports in $\mathcal B$ gives $X$ by the
block-design and sign-change construction in
\cite[Theorem~3.3 and Sections~4.2 and~4.4]{Victoir2004}.
The equal-weight averages of $x_i^2,x_i^4,x_i^2x_j^2$ ($i\ne j$)
are $1/10,1/40,1/120$, respectively, which are the corresponding
moments of normalized surface measure on $\bS^9$.
All monomials having an odd exponent have average zero by sign
symmetry. Hence $X$ is a spherical $5$-design with $240$ points.

Consider the two collections
\[
\begin{array}{c|c}
\mathcal T^+&\mathcal T^-\\\hline
\{1,2,6,8\}&\{1,2,3,6\}\\
\{1,3,7,10\}&\{1,7,8,10\}\\
\{2,3,5,6\}&\{2,5,6,8\}\\
\{5,7,8,10\}&\{3,5,7,10\}
\end{array}
\]
with $\mathcal T^+\subset\mathcal B$ and
$\mathcal T^-\cap\mathcal B=\varnothing$.
The incidence vectors of these blocks are obtained from the two
legs of the OA trade in Example~\ref{exam:binary-oa} by the same
operations: repeating coordinates, taking complements of selected
coordinates, and adjoining zero coordinates.
The identity
\[
 \sum_{B\in\mathcal T^+}\prod_{i\in B}z_i
 -\sum_{B\in\mathcal T^-}\prod_{i\in B}z_i
 =(z_8-z_3)(z_1-z_5)(z_2z_6-z_7z_{10})
\]
shows that they form a $2$-trade: all derivatives of total order
at most two vanish at $z_1=\cdots=z_{10}=1$, so the incidence
counts of every subset of size at most two agree.
This is a minimum-volume $2$-trade of the standard form
described in
\cite[Lemma~4]{GhorbaniKamaliKhosrovshahiKrotov2020},
with four blocks in each leg.
Trades between $2$-$(10,4,2)$ designs were enumerated by
Ramsay~\cite[Section~4.2 and Tables~4--5]{Ramsay1997}.

Put $Y^\pm=\bigcup_{B\in\mathcal T^\pm}Q_B$.
For every monomial of degree at most $5$, its sum over $Q_B$
is either zero or a constant times the indicator that a specified
set of at most two coordinates is contained in $B$.
Thus $Y^+\sqcup Y^-$, with weights $+1$ on $Y^+$ and $-1$ on
$Y^-$, is a spherical null $5$-design.
The same sign-averaging argument applies to any combinatorial
$2$-trade with blocks of size $k$ on $[v]$, using all signed
vectors with those supports and nonzero coordinates $\pm1/\sqrt{k}$.
Assigning each point the signed multiplicity of its block gives
a spherical null $5$-design on $\bS^{v-1}$.

Since $Y^+\subset X$ and $Y^-\cap X=\varnothing$,
Theorem~\ref{thm:switching1}, with the identity embedding,
$T^m=\R^{10}$, $z_0=0$ and $c=1/240$, gives a spherical
$5$-design
\[
 X'=(X\setminus Y^+)\sqcup Y^-
\]
with $240$ points, each of weight $1/240$.
Exactly $64$ points are replaced and $176$ are retained.

Both legs span the coordinate subspace
$V=\operatorname{span}\{e_i:i\in I\}$, where
$I=\{1,2,3,5,6,7,8,10\}$.
Sign symmetry and the fact that each $i\in I$ occurs in two of
the four blocks show that each leg is a spherical $3$-design
on the unit sphere of $V$, identified with $\bS^7$.
Neither leg is a spherical $5$-design on this sphere, since
\[
 \frac1{64}\sum_{y\in Y^\pm}y_i^4
 =\frac1{32}\ne\frac{3}{8(8+2)}\qquad(i\in I),
\]
where the last expression is the corresponding spherical fourth
moment. Their matching moments through degree $5$ nevertheless
allow the complete switching above.

The configurations $X$ and $X'$ are noncongruent.
To see this, write $\mathcal C(\mathcal D)$ for the binary span of the
incidence vectors of a block collection $\mathcal D$, and put
$\mathcal B'=(\mathcal B\setminus\mathcal T^+)\cup\mathcal T^-$.
The $4$-cycles belong to the even-weight subspace of the binary
cycle space of $K_5$. This subspace has dimension $5$, since the
cycle space has dimension $10-5+1=6$ and contains an odd cycle.
The five retained blocks
\[
 \begin{gathered}
 \{1,2,7,9\},\quad \{1,3,5,8\},\quad \{1,4,5,9\},\\
 \{1,4,6,10\},\quad \{2,3,9,10\}
 \end{gathered}
\]
are independent over $\mathbb F_2$: their incidence submatrix
on coordinates $1,2,3,5,6$ has determinant $1$ over $\mathbb F_2$.
Hence $\mathcal C(\mathcal B\setminus\mathcal T^+)=\mathcal C(\mathcal B)$ has
dimension $5$. Each new block is obtained from the corresponding
removed block in the table by interchanging coordinates $3$ and
$8$, so their incidence vectors differ by $e_3+e_8$ over
$\mathbb F_2$. This vector is not in the cycle space, since the
subgraph with edges $ad$ and $cd$ has odd degree at $a$ and $c$.
Therefore $\mathcal C(\mathcal B')=\mathcal C(\mathcal B)+\langle e_3+e_8\rangle$
has dimension $6$.

For $Z=X$ or $X'$, let $L_Z=\operatorname{span}_{\Z}(2Z)$.
Every coordinate occurs in a block, so differences between sign
choices generate every $2e_i$. Thus $L_Z$ is the inverse image
of its block-incidence code under reduction modulo $2$, and
\[
 \operatorname{covol}(L_X)=2^{10-5}=32,
 \qquad \operatorname{covol}(L_{X'})=2^{10-6}=16.
\]
Both configurations have centroid zero, so a congruence between
them would be orthogonal. Such a transformation would take
$L_X$ to $L_{X'}$ and preserve covolume, which is impossible.
\end{exam}

\section{Concluding remarks}
\label{sect:6}
The preceding results give an elementary description of minimum-size
spherical null configurations, including their weights.
The explicit weights yield switching on circles under suitable
placement and weight conditions, while combinatorial trades supply
further examples, including an equal-weight replacement between
noncongruent spherical $5$-designs.

From a different perspective, the null-design framework also provides
a natural interpretation of the Prouhet--Tarry--Escott (PTE) problem.
Let
\[
 A=\{\boldsymbol{a}_1,\ldots,\boldsymbol{a}_m\},
 \qquad
 B=\{\boldsymbol{b}_1,\ldots,\boldsymbol{b}_m\}
 \subset\mathbb Z^r
\]
be a nontrivial solution of the $r$-dimensional PTE problem of degree
$t$~\cite{AlpersTijdeman2007}.
Assume that neither side has repeated points and that common points
have been cancelled, so that $A\cap B=\varnothing$.
Then the sums of all monomials of total degree at most $t$ over
$A$ and $B$ agree.
Set $X=A\sqcup B$ and define $w:X\to\{-1,1\}$ by
\[
 w(x)=
 \begin{cases}
  1,&x\in A,\\
  -1,&x\in B.
 \end{cases}
\]
Then
\[
 \sum_{x\in X}w(x)f(x)=0
 \qquad
 \text{for every }f\in\R_t[x_1,\ldots,x_r].
\]
Thus $(X,w)$ is a Euclidean null $t$-design with integral points
and weights in $\{-1,1\}$.

This interpretation also distinguishes the linear-algebraic existence
of null designs from the arithmetic problem of constructing PTE
solutions with restricted multiplicities.
For rational point configurations with the usual polynomial structure,
the evaluation matrix has rational entries, so any nontrivial
kernel contains a nonzero vector that can be scaled to have
integral entries.
Requiring all nonzero weights to belong to $\{-1,1\}$, however,
imposes an additional restriction not ensured by
Proposition~\ref{prop:existence}.

This viewpoint is complementary to recent connections between PTE
problems and geometric or combinatorial designs.
In particular, PTE solutions arising from ellipsoidal designs and
from combinatorial designs
(see, e.g.,~\cite{Inagaki2026,MatsumuraSawa2025})
can be regarded as structured examples of the general null-design
phenomenon considered here.

It remains to determine which nonempty proper subsets of a given
spherical design can be replaced while respecting its weights. For equal-weight switching of
$Z\subset \bS^{n-1}$, this asks for $A\subsetneq Z$ and a nonempty
configuration $B\subset \bS^{n-1}\setminus Z$ such that
$\sum_{a\in A}f(a)=\sum_{b\in B}f(b)$ for every $f\in P_t(\R^n)$.
Then $(Z\setminus A)\cup B$ has the same moments and cardinality.
The replacement may introduce new points; only $A$ must belong
to the original design. Real null weights on a subset, guaranteed
for sufficiently many points by Corollary~\ref{cor:spherical-null-existence},
do not by themselves provide an equal-weight replacement.

To identify small exchange structures, one may study minimum
null supports with a prescribed affine dimension.
For $t\ge1$, the smallest spherical null $t$-designs have affine
dimension $2$; with equal absolute weights they consist of two
regular $(t+1)$-gons, as illustrated by the cube example.
For higher dimensions, one may seek the minimum size and the
corresponding configurations and weights subject to
\[
 X\subset \bS^{n-1},\qquad \dim\operatorname{aff}X=d,\qquad 3\le d\le n.
\]
Here $X$ is the combined null support of an exchange.
Arbitrary nonzero real weights and equal absolute weights give
minimization problems that should be studied separately.
Applying the resulting configurations to a given design would
still require compatible placement and weights.

\section*{Acknowledgments}
This work was partially supported by JSPS KAKENHI
Grant Numbers 22K03402, 24K06688, 24K06871, 25K17298, 26K00617,
and JST SPRING Grant Number JPMJSP2148.

 \bigskip

 \noindent 
 \textbf{Data Availability Statement.} 
This article is purely theoretical and does not involve any datasets. 

\bigskip

\noindent
\textbf{Conflict of Interest.} 
The authors state that there is no conflict of interest. 

\bigskip

\noindent
\textbf{Declaration of Generative AI Use.}
Generative AI tools were used to assist with mathematical discussions
and manuscript preparation.
The authors take full responsibility for the mathematical content
and accuracy of the paper.

\bibliographystyle{plain}
\bibliography{null_design_v2}

\end{document}